\documentclass[12pt]{article}
\usepackage{graphicx}
\usepackage{amsmath,amsthm,amssymb,enumerate}
\usepackage{euscript,mathrsfs}
\usepackage[left=2cm,right=2cm,top=3.5cm,bottom=3.5cm]{geometry}
\usepackage{color}

\usepackage{soul}

\catcode`\@=11 \@addtoreset{equation}{section}

\catcode`\@=12

\allowdisplaybreaks

\newtheorem{Theorem}{Theorem}[section]
\newtheorem{Proposition}[Theorem]{Proposition}
\newtheorem{Lemma}[Theorem]{Lemma}
\newtheorem{Corollary}[Theorem]{Corollary}

\theoremstyle{definition}
\newtheorem{Definition}[Theorem]{Definition}

\newtheorem{Remark}[Theorem]{Remark}

\newcommand{\bTheorem}[1]{
\begin{Theorem} \label{T#1} }
\newcommand{\eT}{\end{Theorem}}

\newcommand{\bProposition}[1]{
\begin{Proposition} \label{P#1}}
\newcommand{\eP}{\end{Proposition}}

\newcommand{\bLemma}[1]{
\begin{Lemma} \label{L#1} }
\newcommand{\eL}{\end{Lemma}}

\newcommand{\bCorollary}[1]{
\begin{Corollary} \label{C#1} }
\newcommand{\eC}{\end{Corollary}}

\newcommand{\bRemark}[1]{
\begin{Remark} \label{R#1} }
\newcommand{\eR}{\end{Remark}}

\newcommand{\bDefinition}[1]{
\begin{Definition} \label{D#1} }
\newcommand{\eD}{\end{Definition}}

\newcommand{\vrE}{\vr_E}

\newcommand{\Ds}{\mathbb{D}_x}

\newcommand{\vuE}{\vc{u}_E}

\newcommand{\bfphi}{\boldsymbol{\varphi}}

\newcommand{\bFormula}[1]{
\begin{equation} \label{#1}}
\newcommand{\eF}{\end{equation}}

\newcommand{\Ov}[1]{\overline{#1}}

\newcommand{\vr}{\varrho}

\newcommand{\vu}{\vc{u}}
\newcommand{\vm}{\vc{m}}

\newcommand{\vn}{\vc{n}}

\newcommand{\vc}[1]{{\bf #1}}

\newcommand{\Div}{{\rm div}_x}
\newcommand{\Grad}{\nabla_x}

\newcommand{\dx}{\,{\rm d} {x}}

\newcommand{\dt}{\,{\rm d} t }

\newcommand{\intO}[1]{\int_{\Omega} #1 \ \dx}

\newcommand{\D}{{\rm d}}

\newcommand{\ep}{\varepsilon}

\def\softd{{\leavevmode\setbox1=\hbox{d}%
          \hbox to 1.05\wd1{d\kern-0.4ex{\char039}\hss}}}
\definecolor{Cgrey}{rgb}{0.85,0.85,0.85}
\definecolor{Cblue}{rgb}{0.50,0.85,0.85}
\definecolor{Cred}{rgb}{1,0,0}
\definecolor{fancy}{rgb}{0.10,0.85,0.10}

\newcommand\Cbox[2]{%
    \newbox\contentbox%
    \newbox\bkgdbox%
    \setbox\contentbox\hbox to \hsize{%
        \vtop{
            \kern\columnsep
            \hbox to \hsize{%
                \kern\columnsep%
                \advance\hsize by -2\columnsep%
                \setlength{\textwidth}{\hsize}%
                \vbox{
                    \parskip=\baselineskip
                    \parindent=0bp
                    #2
                }%
                \kern\columnsep%
            }%
            \kern\columnsep%
        }%
    }%
    \setbox\bkgdbox\vbox{
        \color{#1}
        \hrule width  \wd\contentbox %
               height \ht\contentbox %
               depth  \dp\contentbox
        \color{black}
    }%
    \wd\bkgdbox=0bp%
    \vbox{\hbox to \hsize{\box\bkgdbox\box\contentbox}}%
    \vskip\baselineskip%
}

\date{}

\begin{document}


\title{On convergence to equilibria of flows of compressible viscous fluids under in/out--flux boundary conditions}

\author{Jan B\v rezina \and Eduard Feireisl
\thanks{The work of E.F. was partially supported by the
Czech Sciences Foundation (GA\v CR), Grant Agreement
18--05974S. The Institute of Mathematics of the Academy of Sciences of
the Czech Republic is supported by RVO:67985840.} 
\and Anton\' \i n Novotn\' y { \thanks{The work of A.N. was supported by Brain Pool program funded by the Ministry of Science and ICT through the National Research Foundation of Korea (NRF-2019H1D3A2A01101128).}}
}


\maketitle

\centerline{Faculty of Arts and Science, Kyushu University;}
\centerline{744 Motooka, Nishi-ku, Fukuoka, 819-0395, Japan}
\centerline{brezina@artsci.kyushu-u.ac.jp}

\centerline{and}

\centerline{Institute of Mathematics of the Academy of Sciences of the Czech Republic;}
\centerline{\v Zitn\' a 25, CZ-115 67 Praha 1, Czech Republic}

\centerline{Institute of Mathematics, Technische Universit\"{a}t Berlin,}
\centerline{Stra{\ss}e des 17. Juni 136, 10623 Berlin, Germany}
\centerline{feireisl@math.cas.cz}

\centerline{and}

\centerline{IMATH, EA 2134, Universit\'e de Toulon,}
\centerline{BP 20132, 83957 La Garde, France}
\centerline{novotny@univ-tln.fr}

\begin{abstract}

We consider the barotropic Navier--Stokes system describing the motion of a compressible Newtonian fluid in a bounded domain 
with in and out flux boundary conditions. We show that if the boundary velocity coincides with that of a rigid motion, all solutions 
converge to an equilibrium state for large times.

\end{abstract}

{\bf Keywords:} compressible Newtonian fluid, Navier--Stokes system, in/out--flux boundary conditions, long--time behavior

{\bf MSC:} 
\bigskip

\section{Introduction}
\label{i}

The barotropic Navier--Stokes system: 
\begin{equation} \label{i1}
\begin{split}
\partial_t \vr + \Div (\vr \vu) &= 0,\\
\partial_t (\vr \vu) + \Div (\vr \vu \otimes \vu) + \Grad p(\vr) &= 
\Div \mathbb{S}(\Ds \vu) + \vr \Grad G,\\ 
\mathbb{S}(\Ds \vu) &= \mu \left( \Grad \vu + \Grad^t \vu - \frac{2}{d} \Div \vu \mathbb{I} \right) + 
\lambda \Div \vu \mathbb{I},\ \mu > 0, \ \lambda \geq 0,\\ \mbox{with}\ \Ds \vu &\equiv \frac{1}{2} \Big( \Grad \vu + \Grad^t \vu \Big), 
\end{split}
\end{equation}
is a well--established model in continuum fluid mechanics governing the time evolution of the mass density $\vr = \vr(t,x)$ and the 
velocity $\vu = \vu(t,x)$ of a compressible viscous fluid. In the fluid if confined to a bounded domain $\Omega \subset R^d$, $d=1,2,3$, 
suitable boundary conditions must be prescribed to obtain a well posed problem. Here we consider the realistic situation with a given 
boundary velocity, 
\begin{equation} \label{i2}
\vu|_{\partial \Omega} = \vu_b, 
\end{equation}
and, decomposing the boundary as
\[
\partial \Omega = \Gamma_{\rm in} \cup \Gamma_{\rm out}, 
\ \Gamma_{\rm in} = \left\{ x \in \partial \Omega \ \Big|\
\ \mbox{the outer normal}\ \vc{n}(x) \ \mbox{exists, and}\ \vu_B(x) \cdot \vc{n}(x) < 0 \right\}, 
\]
we prescribe the density on the in--flow component, 
\begin{equation} \label{i3}
\vr|_{\Gamma_{\rm in}} = \vr_b.
\end{equation}
Our goal is to describe the long--time behavior of finite energy weak solutions to the problem \eqref{i1}--\eqref{i3}.

Note that the long--time behavior of solutions is well understood under the no--slip boundary conditions $\vu_b \equiv 0$, see 
\cite{FP7}, \cite{FP9}, \cite{NOS1}, \cite{NOST} for general results if $d=2,3$ and Melinand and Zumbrun \cite{MelZum} for refined arguments if $d=1$. The $\omega$--limit set of 
any solution trajectory $t \mapsto [\vr(t, \cdot), (\vr \vu)(t, \cdot)]$ is contained in the set of stationary (static) solutions
$[\vrE, 0]$, 
\begin{equation} \label{i4}
\Grad p(\vrE) = \vr_E \Grad G,\ \vr_E \geq 0,\ \intO{\vrE} = \intO{\vr(0, \cdot)} = M_0.
\end{equation}
If the problem \eqref{i4} admits a unique solution, any trajectory converges to it. The same is true if the set of solutions 
of \eqref{i4} consists of isolated points. The case when \eqref{i4} admits a continuum of solutions  
remains an outstanding open problem. Note that in this case the equilibria $\vrE$ necessarily contain vacuum, meaning $\vrE$ vanishes on a set of non--zero measure, see \cite{FP7}. 

Much less is known in the case of non--trivial in/out flow velocity. Melinand and Zumbrun \cite{MelZum} studied the problem in the mono--dimensional 
case $d=1$ and with $G = 0$ in the framework of strong solutions. They show (non--linear) stability of the stationary solutions with constant velocity $\vu_b$ and their small perturbations. They also show that linear stability implies nonlinear stability in the 
general case. 

Motivated by \cite{FP9}, we study stability and convergence to the static states in the multi--dimensional case, with the velocity $\vuE$ associated to a 
\emph{rigid motion}, meaning 
\begin{equation} \label{i5}
\Ds \vuE  = 0.
\end{equation}
The corresponding density $\vrE$ satisfies 
\begin{equation} \label{i6}
\begin{split}
\Div (\vrE \vuE) &= 0,\\ 
\Div (\vrE \vuE \otimes \vuE) + \Grad p(\vrE) &= \vrE \Grad G.
\end{split}
\end{equation}
Accordingly, we consider the problem \eqref{i1}--\eqref{i3} with the boundary conditions 
\begin{equation} \label{i7}
\vu_b = \vuE, \ \vr_b = \vrE.
\end{equation} 

Under the hypothesis \eqref{i7}, and if the stationary density $\vrE$ is strictly positive, the problem \eqref{i1}--\eqref{i3} admits a Lyapunov function, namely the relative energy
\[
\intO{ E\left(\vr, \vu \Big| \vrE, \vuE \right) },\   E\left(\vr, \vu \Big| \vrE, \vuE \right) 
\equiv \left[ \frac{1}{2} \vr |\vu - \vuE|^2 + P(\vr) - P'(\vrE)(\vr - \vrE) - P(\vrE)  \right],
\]
see Section \ref{S}. The situation becomes more delicate if $\vrE$ vanishes on a non--trivial part of $\Omega$. In that case, the stationary problem may admit more (infinitely many) solutions even if the total mass is prescribed.  

Our main result asserts that any \emph{weak} solution of the problem \eqref{i1}--\eqref{i3}, satisfying 
a suitable form of energy inequality, approaches the equilibrium solution 
$[\vrE, \vuE]$ as $t \to \infty$ as long as the stationary problem \eqref{i6} admits a unique solution. To the best of our knowledge, this is the first result of this kind in the multi--dimensional case under the non--zero in/out flow boundary conditions. Note 
that such a result does not follow from ``standard'' arguments, even if $\vrE > 0$, as the Lyapunov function
\[ 
t \mapsto \intO{ E\left(\vr, \vu \Big| \vrE, \vuE \right)(t, \cdot) } 
\]
is not continuous on the trajectories generated by weak solutions. In addition, we show that the convergence is uniform with respect to bounded energy initial data. 

The paper is organized as follows. In Section \ref{M}, we recall the concept of weak solution to the Navier--Stokes system and state our main result. Section \ref{S} is devoted to the stationary problem \eqref{i6}. In particular, we establish several conditions sufficient for its unique solvability. The main convergence result is shown in Section \ref{C}. 

\section{Weak solutions, energy inequality, main results}
\label{M}

We start by introducing the main hypotheses imposed on the structural properties of the potential $G$ and the pressure $p$.
In what follows, we shall always assume that $\Omega \subset R^d$ is a bounded Lipschitz domain.
Keeping in mind the iconic example of the gravitational potential, we require only 
\begin{equation} \label{MH1}
G \in C^1(\Ov{\Omega}).
\end{equation} 
As for the pressure, we assume 
\begin{equation} \label{MH2}
p \in C^1[0, \infty),\ p(0) = 0, \ p'(\vr) > 0 \ \mbox{for}\ \vr > 0, \ p'(\vr) \approx \vr^{\gamma - 1}, \ \gamma > 1 
\ \mbox{as}\ \vr \to \infty.
\end{equation}
Here, the symbol $p'(\vr) \approx \vr^{\gamma - 1}$ as $\vr \to \infty$ means 
\[
\underline{p} \vr^{\gamma - 1} \leq p'(\vr) \leq \Ov{p} \vr^{\gamma - 1}  \ \mbox{for all}\ 
\vr > 1, \ \mbox{where}\ \underline{p} > 0.
\]
Accordingly, the pressure potential $P$ defined as 
\[
P'(\vr) \vr - P(\vr) = p(\vr),\ P(0) = 0, \ \Rightarrow \ P''(\vr) = \frac{p'(\vr)}{\vr} \ \mbox{for}\ \vr > 0,
\]
is a strictly convex function on $[0, \infty)$. Without loss of generality, we may therefore assume 
\[
P'(\vr) \to - \infty \ \mbox{if} \ \vr \to 0+ \ \mbox{or}\ P'(\vr) \to 0 \ \mbox{if}\ \vr \to 0+, 
\]
adding a linear function to $P$ in the latter case if necessary.

\subsection{Weak solutions to the Navier--Stokes system}
\label{MS1}

The functions $[\vr, \vu]$ represent a weak solution of the Navier--Stokes system \eqref{i1}--\eqref{i3} in $[0, \infty) \times \Omega$, with the boundary data 
\[
\vu_b = \vuE|_{\partial \Omega},\ \vr_b = \vrE|_{\partial \Omega},
\]
if: 
\begin{itemize}
\item
\[
\begin{split}
\vr \in C_{{\rm weak,loc}}([0, \infty); L^\gamma (\Omega)),\ \vr \geq 0,\  
\vm \equiv \vr \vu &\in C_{{\rm weak,loc}}([0, \infty); L^{\frac{2 \gamma}{\gamma + 1}}(\Omega; R^d)),\\ 
(\vu - \vuE) \in L^2_{\rm loc}([0,\infty); W^{1,2}(\Omega; R^d)),\ \vr &\in L^\gamma_{{\rm loc}}([0, \infty); 
L^\gamma (\Gamma_{\rm out}; \D |\vu_b \cdot \vc{n}|)).
\end{split}
\]

\item
Equation of continuity
\begin{equation} \label{M1}
\begin{split}
\left[ \intO{ \vr \varphi } \right]_{t = 0}^{t = \tau} &+ 
\int_0^\tau \int_{\Gamma_{\rm out}} \varphi \vr \vuE \cdot \vc{n} \ \D \ S_x 
+ 
\int_0^\tau \int_{\Gamma_{\rm in}} \varphi \vrE \vuE \cdot \vc{n} \ \D \ S_x\\ &= 
\int_0^\tau \intO{ \Big[ \vr \partial_t \varphi + \vr \vu \cdot \Grad \varphi \Big] } \dt 
\end{split}
\end{equation}
holds for any $0 \leq \tau < \infty$, and any test function
for any $\varphi \in C^1_c([0,\infty) \times \Ov{\Omega})$. 

In addition, we require also the renormalized version of \eqref{M1},
\begin{equation} \label{M1a}
\begin{split}
\left[ \intO{ b(\vr) \varphi } \right]_{t = 0}^{t = \tau} &+ 
\int_0^\tau \int_{\Gamma_{\rm out}} \varphi b(\vr) \vuE \cdot \vc{n} \ \D \ S_x 
+ 
\int_0^\tau \int_{\Gamma_{\rm in}} \varphi b(\vrE) \vuE \cdot \vc{n} \ \D \ S_x\\ &= 
\int_0^\tau \intO{ \Big[ b(\vr) \partial_t \varphi + b(\vr) \vu \cdot \Grad \varphi - 
\Big( b'(\vr) \vr - b(\vr) \Big) \Div \vu \Big] } \dt 
\end{split}
\end{equation}
to be satisfied for any $0 \leq \tau < \infty$, any test function
for any $\varphi \in C^1_c([0,\infty) \times \Ov{\Omega})$, and any $b \in C^1[0, \infty)$, $b' \in C_c[0, \infty)$.

\item 
Momentum equation
\begin{equation} \label{M2}
\begin{split}
\left[ \intO{ \vr \vu \cdot \bfphi } \right]_{t=0}^{t = \tau} &= 
\int_0^\tau \intO{ \Big[ \vr \vu \cdot \partial_t \bfphi + \vr \vu \otimes \vu : \Grad \bfphi 
+ p(\vr) \Div \bfphi - \mathbb{S}(\Ds \vu) : \Grad \bfphi \Big] }\\
&+ \int_0^\tau \intO{\vr \Grad G \cdot \bfphi    } \dt
\end{split}
\end{equation}
holds for any $0 \leq \tau < \infty$, and any test function
$\bfphi \in C^1_c([0,\infty) \times {\Omega}; R^d)$.

\end{itemize}

\subsection{Energy balance} 

The energy inequality is an indispensable part of the definition of weak solution.
In view of direct calculations presented in the Appendix it takes the form 
\begin{equation} \label{M3}
\begin{split}
-&\int_0^\infty \partial_t \psi  \intO{\left[ \frac{1}{2} \vr |\vu - \vuE|^2 + P(\vr) \right] }\dt  + 
\int_0^\infty \psi \intO{ \mathbb{S}(\Ds \vu) : \Ds \vu } \dt \\  
&+ \int_0^\infty \psi \int_{\Gamma_{\rm out}} P(\vr)  \vuE \cdot \vc{n} \ \D S_x \dt  +
\int_0^\infty \psi \int_{\Gamma_{\rm in}} P(\vr_E)  \vuE \cdot \vc{n} \ \D S_x \dt 
\\	
&\leq
\psi(0) \intO{\left[ \frac{1}{2} \vr(0, \cdot)  | \vu(0, \cdot) - \vuE|^2  + P(\vr(0, \cdot)) \right]} \\
&- 
\int_0^\infty \psi \intO{ \left[ \vr \vu \otimes \vu + p(\vr) \mathbb{I} \right]  :  \Grad \vuE } \dt + \int_0^\infty \psi \intO{ {\vr} { \vu \cdot \frac12 \Grad |\vuE|^2}   } 
\dt\\ &+ \int_0^\infty \psi \intO{ \mathbb{S}(\Ds \vu) : \Ds \vuE } \dt + \int_0^\infty \psi \intO{ \vr \Grad G \cdot (\vu - \vu_E) }\dt 
\end{split}
\end{equation}
for any $\psi \in C^1_c[0, \infty)$, $\psi \geq 0$.

\begin{Remark} \label{MR3}

The energy can be defined in terms of the density and momentum that are weakly continuous quantities in time:
\[
E \left(\vr, \vu \ \Big| \vuE \right) \equiv \left[ \frac{1}{2} \vr |\vu - \vuE|^2 + P(\vr) \right] = 
\left[ \frac{1}{2} \frac{|\vm|^2}{\vr} - \vm \cdot \vuE + \frac{1}{2} \vr |\vuE|^2 + P(\vr) \right], \ \vm \equiv 
\vr \vu. 
\]
Moreover, with the convention 
\[
E \left(\vr, \vu \ \Big| \vuE \right) = \infty \ \mbox{if}\ \vr < 0 \ \mbox{or}\ \vr = 0, \vm \ne 0,\ 
E \left(\vr, \vu \ \Big| \vuE \right) = 0 \ \mbox{if}\ \vr = 0, \ \vm = 0,
\]
$E$ is a convex l.s.c. function of $[\vr, \vm] \in R^{d + 1}$.

\end{Remark}

\begin{Definition}[Finite energy weak solution] \label{MD1}

A weak solution $[\vr, \vu]$ specified in Section \ref{MS1} satisfying the energy inequality 
\eqref{M3} is called \emph{finite energy weak solution} of the Navier--Stokes system \eqref{i1}--\eqref{i3} in 
$[0, \infty) \times \Omega$.

\end{Definition}

The existence of finite energy weak solutions for the Navier--Stokes system with in/out flux boundary conditions has been proved 
in \cite{ChJiNo}, \cite{ChNoYa}, \cite{KwoNovSat} (see also Girinon \cite{GI}) under additional assumptions on smoothness of the domain $\Omega$ and for $\gamma > \frac{d}{2}$.
At this stage the total energy 
\[
\intO{ E \left( \vr, \vu \ \Big| \vuE \right) }
\]
is not necessarily a decreasing function of time. Further assumptions on $\vuE$ and $\vr_b$ specified below are necessary to convert 
it to a kind of Lyapunov function for the system.

\subsection{Main result}

We are ready to state our main result. 

\begin{Theorem}[Convergence to equilibrium] \label{MT1}

Let $\Omega \subset R^d$, $d=2,3$ be a bounded Lipschitz domain. Let $G$ and $p$ satisfy the hypotheses \eqref{MH1}, \eqref{MH2}, with 
\[
\gamma > \frac{d}{2}.
\]
Let $\vuE$ be a given field such that 
\begin{equation} \label{M4}
\Ds \vuE = 0,\ \Grad G \cdot \vuE = 0.
\end{equation}
Let $\vrE$ be a density field solving the stationary problem \eqref{i6} with the given $\vuE$ such that 
\[
\vrE \geq 0, \ \mbox{the set}\ \left\{ x \in \Omega \ \Big| \ \vrE(x) > 0 \right\} \ne \emptyset 
\ \mbox{is connected in}\ \Omega, \ \vrE|_{\Gamma_{\rm in}} > 0.
\] 

Let $[\vr, \vu]$ be a finite energy weak solution of the problem \eqref{i1}--\eqref{i3} in $[0, \infty) \times \Omega$, with 
the boundary conditions \eqref{i7}, and 
\[
\intO{ E \left( \vr, \vu \ \Big|\ \vuE \right)(0, \cdot) } \leq E_0, 
\]
\[
\intO{ \vr(0, \cdot) } = M_0 > 0, \ M_0 = \intO{\vrE} \ \mbox{if}\ \Gamma_{\rm in} = \emptyset.
\]

Then for any $\ep > 0$, there exists $T = T(\ep)$ depending only on $E_0$ such that 
\[
\| \vr(t, \cdot) - \vrE \|_{L^\gamma(\Omega)} + 
\| \vr (\vu - \vuE)(t, \cdot) \|_{L^{\frac{2 \gamma}{\gamma + 1}}(\Omega; R^d)} < \ep 
\ \mbox{for all}\ t > T(\ep).
\]
 
\end{Theorem}

\begin{Remark} \label{MR1}

It follows from the equation \eqref{i6} that the pressure $p(\vrE)$ is a continuously differentiable function in $\Ov{\Omega}$, in particular, $\vrE \in C(\Ov{\Omega})$.

\end{Remark}

\begin{Remark} \label{MR2}

As $\Div \vuE = 0$, we have 
\[
\int_{\partial \Omega} \vuE \cdot \vc{n}\ \D S_x = 0.
\]
Consequently, if $\Gamma_{\rm in} = \emptyset$, then necessarily 
\begin{equation} \label{M12}
\vuE \cdot \vc{n}|_{\partial \Omega} = 0.
\end{equation}
As $\vuE$ is the velocity of a rigid motion and $\Omega$ is bounded, relation \eqref{M12} implies 
either $\Omega$ is rotationally symmetric or $\vuE = 0$. In both cases, the total mass 
\[
\intO{ \vr(t, \cdot) } = M_0 \ \mbox{is a constant of motion.}
\]

\end{Remark}

The following two sections are devoted to the proof of Theorem \ref{MT1}. 

\section{Stationary problem}
\label{S}

The energy inequality \eqref{M3} simplifies to
\begin{equation} \label{M9}
\begin{split}
&-\int_0^\infty \partial_t \psi \intO{\left[ \frac{1}{2} \vr |\vu - \vuE|^2 + P(\vr) - \vr \left( \frac{1}{2} |\vuE|^2 + G \right) \right] } \dt + 
\int_0^\infty \psi \intO{ \mathbb{S}(\Ds \vu) : \Ds \vu } \dt\\  
&+\int_0^\infty \psi \int_{\Gamma_{\rm out}} \left[ P(\vr) - (\vr - \vrE) \left(\frac{1}{2} |\vuE|^2 + G \right) - 
P(\vrE) \right] \vuE \cdot \vc{n} \ \D S_x \dt\\&+
\int_0^\infty \psi \int_{\partial \Omega} \left[ P(\vr_E) - \vrE \left( \frac{1}{2} |\vuE|^2 + G \right) \right]  \vuE \cdot \vc{n} \ \D S_x \dt \\
&
\leq \psi(0) \intO{ \left[ E \left(\vr, \vu \ \Big|\ \vuE \right) +  \left( \frac{1}{2} |\vuE|^2 + G \right)\right](0, \cdot) }
\end{split}
\end{equation}
for any $\psi \in C^1_c[0, \infty)$, $\psi \geq 0$.

\subsection{Stationary equation of continuity}

Next we use the hypothesis that the boundary data for the density $\vr_b = \vrE|_{\Gamma_{\rm in}}$ are determined by the stationary 
density $\vrE$ satisfying, in particular, the equation of continuity
\begin{equation} \label{M10}
\Div(\vrE \vuE) = 0 \ \mbox{in}\ \mathcal{D}'(\Omega).
\end{equation}

It follows from \eqref{M10} that 
\[
\int_{\partial \Omega}  \vrE \left( \frac{1}{2} |\vuE|^2 + G \right)    \vuE \cdot \vc{n} \ \D S_x 
= - \intO{ \vrE \Grad  \left( \frac{1}{2} |\vuE|^2 + G \right)  \cdot \vuE } = 0,
\]
where the last equality follows from \eqref{M4} and 
\begin{equation} \label{M8}
\Grad |\vuE|^2 \cdot \vuE = - 2 \vuE \cdot \Grad \vuE \cdot = 
- \vuE \cdot \Grad \vuE \cdot \vuE - \vuE \cdot \Grad^t \vuE \cdot \vuE = - 2 \vuE \cdot \Ds \vuE \cdot \vuE  = 0.
\end{equation} 
Similarly, using $\Div \vuE = 0$ we get by renormalization 
\[
\Div (P(\vrE) \vuE ) = 0 \ \Rightarrow \ \int_{\partial \Omega} P(\vrE) \vuE \cdot \vc{n}\ \D S_x = 0.
\]
Consequently, the energy inequality \eqref{M9} takes the form 
\begin{equation} \label{M11}
\begin{split}
&- \int_0^\infty \partial_t \psi \intO{\left[ \frac{1}{2} \vr |\vu - \vuE|^2 + P(\vr) - \vr \left( \frac{1}{2} |\vuE|^2 + G \right) \right] } \dt + 
\int_0^\infty \psi \intO{ \mathbb{S}(\Ds \vu) : \Ds \vu } \dt\\  
&+\int_0^\infty \psi \int_{\Gamma_{\rm out}} \left[ P(\vr) - (\vr - \vrE) \left(\frac{1}{2} |\vuE|^2 + G \right) - 
P(\vrE) \right] \vuE \cdot \vc{n} \ \D S_x \dt \\ &\leq
\psi(0) \intO{ \left[ E \left( \vr, \vu \ \Big| \vuE \right) - \vr \left( \frac{1}{2} |\vuE|^2 + G \right) \right] (0, \cdot) } 
\end{split}
\end{equation}
for any $\psi \in C^1_c[0, \infty)$, $\psi \geq 0$.
Note that the result holds under general assumption on $\vrE$, in particular, it is enough that $\vrE \in C(\Ov{\Omega})$,
$\vrE \geq 0$, not necessarily $\vrE > 0$.

\subsection{Stationary momentum equation}
\label{SME}

In view of 
\begin{equation} \label{M7}
\vuE \cdot \Grad \vuE = 2 \vuE \cdot \Ds \vuE - \vuE \cdot \Grad^t \vuE = - \frac{1}{2} \Grad |\vuE|^2,
\end{equation}
 the stationary momentum equation can be written in the form 
\begin{equation} \label{S2}
\Grad p(\vrE) = \vrE \Grad \left( G + \frac{1}{2} |\vuE|^2 \right),
\end{equation}
in particular $p(\vrE) \in C^1(\Ov{\Omega})$, and $\vrE \in C(\Ov{\Omega})$. We point out that $\vrE$ need not be continuously differentiable on the boundary of its domain of positivity. 

If $\vrE > 0$, we can rewrite \eqref{S2} as 
\begin{equation} \label{S3}
\Grad P'(\vrE) = \Grad \left( G + \frac{1}{2} |\vuE|^2 \right) 
\ \Rightarrow \ P'(\vrE) = G + \frac{1}{2} |\vuE|^2 - C_E,
\end{equation} 
where $C_E$ is a constant. In accordance with the hypotheses of Theorem \ref{MT1}, the domain of positivity of $\vrE$, 
\[
\left\{ x \in \Omega \ \Big| \ \vrE(x) > 0 \right\} 
\]
is bounded and connected in $\Omega$; whence $\vrE$ is given through formula 
\begin{equation} \label{S4}
\begin{split}
\vrE(x) &= (P')^{-1} \left[ G(x) + \frac{1}{2} |\vuE(x)|^2 - C_E \right]^+ \ \mbox{if}\ P'(0+) = 0, \\ 
\vrE(x) &= (P')^{-1} \left( G(x) + \frac{1}{2} |\vuE(x)|^2 - C_E \right)\ \mbox{if}\ P'(0+) = -\infty.  
\end{split}
\end{equation}
Note that in the latter case vacuum does not occur, $\vrE > 0$ in $\Omega$. The constant $C_E \in R$ is uniquely determined by the 
boundary value $\vr_b = \vrE|_{\Gamma_{\rm in}}$ if $\Gamma_{\rm in} \ne \emptyset$ or by the total mass 
\[
M_0 = \intO{ \vrE }
\]
in the case $\Gamma_{\rm in} = \emptyset$.

Finally, we rewrite the energy inequality \eqref{M11} in the form 
\begin{equation} \label{S5}
\begin{split}
&- \int_0^\infty \partial_t \psi \intO{\left[ \frac{1}{2} \vr |\vu - \vuE|^2 + P(\vr) - (\vr - 
\vrE) \left( G+ \frac{1}{2} |\vuE|^2  - C_E\right) - P(\vrE) \right] } \dt\\ &+ 
\int_0^\infty \psi \intO{ \mathbb{S}(\Ds \vu) : \Ds \vu } \dt\\  
&+\int_0^\infty \psi \int_{\Gamma_{\rm out}} \left[ P(\vr) - (\vr - \vrE) \left(G+ \frac{1}{2} |\vuE|^2  - C_E\right) - 
P(\vrE) \right] \vuE \cdot \vc{n} \ \D S_x \dt \\ &\leq
\psi(0) \intO{ \left[ E \left( \vr, \vu \ \Big| \vuE \right) - (\vr - \vrE) \left( G + \frac{1}{2} |\vuE|^2  - C_E \right) 
- P(\vrE) \right] (0, \cdot) } 
\end{split}
\end{equation}
for any $\psi \in C^1_c[0, \infty)$, $\psi \geq 0$. Here, we have 
\begin{equation} \label{S6}
\begin{split}
P(\vr) &- (\vr - 
\vrE) \left( G+ \frac{1}{2} |\vuE|^2  - C_E\right) - P(\vrE) \\ &= 
P(\vr) - (\vr - 
\vrE) P'(\vrE) - P(\vrE) \geq 0 \ \mbox{whenever}\ \vrE > 0,
\end{split}
\end{equation}
and 
\begin{equation} \label{S7}
\begin{split}
P(\vr) &- (\vr - 
\vrE) \left( G+ \frac{1}{2} |\vuE|^2  - C_E\right) - P(\vrE) \\ &= 
P(\vr) - \vr \left( G+ \frac{1}{2} |\vuE|^2  - C_E\right)  \geq P(\vr) \geq 0 \ \mbox{if}\ \vrE = 0.
\end{split}
\end{equation}

In particular, the function $\mathcal{E}$,  
\[
\mathcal{E} : t \mapsto \intO{ \left[ \frac{1}{2} \vr |\vu - \vuE|^2 + P(\vr) - (\vr - 
\vrE) \left( G+ \frac{1}{2} |\vuE|^2  - C_E\right) - P(\vrE) \right](t, \cdot) } 
\]
coincides on the set of full measure in $(0, \infty)$ with a non--increasing function and moreover
\begin{equation} \label{S8}
\mathcal{E}(t) \to 0 \ \mbox{as}\ t \to \infty \ \Rightarrow \ 
\| \vr(t, \cdot) - \vrE \|_{L^\gamma(\Omega)} + 
\| \vr (\vu - \vuE)(t, \cdot) \|_{L^{\frac{2 \gamma}{\gamma + 1}}(\Omega; R^d)} \to 0 
\ \mbox{as}\ t \to \infty.
\end{equation}

\section{Convergence to equilibria}
\label{C}

Our goal is to show Theorem \ref{MT1}. We start with the following auxilliary result: 
\begin{Lemma} \label{CL1}
Let $\{ \vr_n, \vu_n \}_{n=1}^\infty$ be a sequence of finite energy weak solutions to the Navier--Stokes system 
\eqref{i1}--\eqref{i3} on a time interval $(0, 1)$ such that 
\[
\intO{ E \left( \vr_n, \vu_n \ \Big| \ \vuE \right) } \leq E_0, 
\]
\[
\int_0^{1} \intO{ \mathbb{S}(\Ds \vu_n) : \Ds \vu_n } \dt \leq E_0 
\]
\mbox{uniformly for}\ n =1,2,\dots, 
\[
\Div \vu_n \to 0 \ \mbox{in}\ L^2(0, 1; L^2(\Omega)).
\]

Then we have
\[
\begin{split}
\vr_n &\to \vr \ \mbox{in}\ L^{\gamma + \alpha}((0,  1) \times \Omega),\\
\vu_n &\to \vu \ \mbox{weakly in}\ L^2(0,1; W^{1,2}(\Omega; R^d)),\\
\vr_n \vu_n \otimes \vu_n &\to \vr \vu \otimes \vu \ \mbox{in}\ L^{1 + \alpha}((0,1) \times \Omega; R^{d \times d})
\end{split}
\]
for a certain $\alpha > 0$, passing to a suitable subsequence as the case may be. 

\end{Lemma}

The proof of Lemma \ref{CL1} is based on nowadays standard arguments of the theory of compressible Navier--Stokes system and may be found in \cite{FeiPr}. 

\subsection{Convergence to equilibria}
\label{CTE}

To show convergence we introduce the sequence of time--shifts:
\[
\vr_n (t,x) = \vr(t + n,x), \ \vu_n(t,x) = \vu(t + n,x),\ n=1,2,\dots 
\]
where $[\vr, \vu]$ is a global--in--time finite energy weak solution to the Navier--Stokes system. It follows from the energy inequality 
\eqref{S5} that 
\[
\int_0^1 \intO{ \mathbb{S}(\Ds \vu_n) : (\Ds \vu_n) } \dt = 
\int_0^1 \intO{ \mathbb{S}(\Ds (\vu_n - \vuE)) : (\Ds (\vu_n - \vuE)) } \dt \to 0 
\ \mbox{as}\ n \to \infty; 
\]
whence, by virtue of Korn--Poincar\' e inequality, 
\[
\vu_n \to \vuE \ \mbox{in}\ L^2(0,1; W^{1,2}(\Omega; R^d)).
\]
Moreover, applying Lemma \ref{CL1} we may perform the limit in the equations \eqref{M1}, \eqref{M2} obtaining 
\begin{equation} \label{C1}
\begin{split}
& 
\int_0^1 \int_{\Gamma_{\rm out}} \varphi \vr \vuE \cdot \vc{n} \ \D \ S_x 
+ 
\int_0^1 \int_{\Gamma_{\rm in}} \varphi \vrE \vuE \cdot \vc{n} \ \D \ S_x\\ &= 
\int_0^1 \intO{ \Big[ \vr \partial_t \varphi + \vr \vuE \cdot \Grad \varphi \Big] } \dt 
\end{split}
\end{equation}
for any test function $\varphi \in C^1_c((0,1) \times \Ov{\Omega})$,
\begin{equation} \label{C2}
-\int_0^1 \intO{ \Big[ \vr \vuE \cdot \partial_t \bfphi + \vr \vuE \otimes \vuE : \Grad \bfphi 
+ p(\vr) \Div \bfphi  \Big] }
= \int_0^1 \intO{\vr \Grad G \cdot \bfphi    } \dt.
\end{equation}
for any test function
$\bfphi \in C^1_c((0,1) \times {\Omega}; R^d)$.
 
It follows from \eqref{C1} that 
\begin{equation} \label{C3}
-\int_0^1 \intO{ \vr \vuE \partial_t \bfphi } \dt = 
\int_0^1 \intO{ \vr \vuE \cdot \Grad \vuE \ \bfphi } \dt + 
\int_0^1 \intO{ \vr \vuE \otimes \vuE : \Grad \bfphi } \dt
\end{equation}
for any $\bfphi \in C^1_c((0,T) \times {\Omega}; R^d)$. In particular, we deduce from   \eqref{C2} using \eqref{M7}  that  
\begin{equation} \label{C1+}
- \int_0^1 \intO{
p(\vr) \Div \bfphi  }
= \int_0^1 \intO{\vr \Grad \left( G + \frac{1}{2} |\vuE|^2 \right) \cdot \bfphi } \dt
\end{equation}
for any $\bfphi \in C^1_c((0,1) \times \Omega; R^d)$. Thus we get 
\[
\Grad p(\vr(t, \cdot)) = \vr(t, \cdot) \Grad \left( G + \frac{1}{2} |\vuE|^2 \right) \ \mbox{in}\ 
\mathcal{D}'(\Omega) \ \mbox{for a.a.}\ t \in (0,1), 
\]
from which, by a simple bootstrap argument, we deduce
\begin{equation} \label{C2+}
p(\vr(t, \cdot)) \in C^1(\Ov{\Omega}),\ 
\vr(t, \cdot) \in C(\Ov{\Omega}), \ \mbox{and}\ \Grad p(\vr(t, \cdot)) = \vr(t, \cdot) \Grad \left( G + \frac{1}{2} |\vuE|^2 \right)
\ \mbox{for a.a.}\ t \in (0,1).
\end{equation}

If $\Gamma_{\rm in} = \emptyset$, then \[
\intO{\vr(t, \cdot)} = M_0
\] for any $t \in (0,1)$ and whence we deduce from {\eqref{C2+}}, exactly as in Section \ref{SME}, that   
\begin{equation} \label{C3+}
\vr(t, \cdot) = \vrE \ \mbox{for a.a. { $t \in (0,1)$}.}
\end{equation}
Similarly, as $\Div \vuE = 0$ and $\vr(t, \cdot)$ is continuous, we deduce from \eqref{C1+} that 
\[
\vr(t, \cdot)|_{\Gamma_{\rm in}} = \vrE \ \mbox{for a.a.}\ t \in (0,1),
\]
which yields the same conclusion \eqref{C3}. 

Consequently, there is a sequence $t_n \to \infty$ such that 
\[
\mathcal{E}(t_n) \to 0 \ \mbox{as}\ t_n \to \infty. 
\]
As $\mathcal{E}(t)$ is non--increasing, this yields the desired conclusion 
\begin{equation} \label{CC1}
\| \vr(t, \cdot) - \vrE \|_{L^\gamma(\Omega)} + 
\| \vr (\vu - \vuE)(t, \cdot) \|_{L^{\frac{2 \gamma}{\gamma + 1}}(\Omega; R^d)} \to 0 
\ \mbox{as}\ t \to \infty.
\end{equation}

\subsection{Uniform convergence}

To show uniform convergence claimed in Theorem \ref{MT1}, it is enough to show 
\[
\mathcal{E}(t) = \intO{\left[ \frac{1}{2} \vr |\vu - \vuE|^2 + P(\vr) - (\vr - 
\vrE) \left( G+ \frac{1}{2} |\vuE|^2  - C_E\right) - P(\vrE) \right](t, \cdot) } \to 0 
\ \mbox{as}\ t \to \infty
\]
uniformly for 
\[
\mathcal{E}(0+) \leq E_0.
\]

Arguing by contradiction, we suppose there is $\delta > 0$, a sequence of time $t_m \to \infty$, and a sequence of global in time solutions 
$\{ \vr_m, \vu_m \}_{m=1}^\infty$, with the associated energies $\mathcal{E}_m$ such that 
\begin{equation} \label{CC2}
\mathcal{E}_m(0+) \leq E_0,\ \mathcal{E}_m(t) \geq \delta > 0 \ \mbox{for any}\ t \in [0, t_m].
\end{equation}
However, as $\mathcal{E}_m$ are non--increasing in time, $\mathcal{E}_m \geq 0$ and satisfying the energy inequality 
\eqref{S5}, we get 
\[
\int_0^{t_m} \intO{ \mathbb{S}(\Ds \vu_m) : \Ds (\vu_m) } \dt \leq E_0 \ \mbox{uniformly for}\ m \to \infty. 
\]
Consequently, there must be another sequence $\tau_m \to \infty$ such that 
\[
(\tau_m, \tau_{m} + 1) \subset [0, t_m),\ \int_{\tau_m}^{\tau_m + 1} \intO{ \mathbb{S}(\Ds \vu_m) : \Ds (\vu_m) } \dt
\to 0 \ \mbox{as}\ m \to \infty.
\]
Thus repeating the arguments of Section \ref{CTE} we would obtain another sequence $s_m \to \infty$, 
$s_m \leq t_m$ such that 
\[
\mathcal{E}(s_m) \to 0 \ \mbox{as}\ m \to \infty 
\]
in contrast with \eqref{CC2}. 

We have proved Theorem \ref{MT1}.

\section{Appendix}

Below we present the formal derivation of energy balance \eqref{M3}.

Since $\vu_E$ is time independent we get from the balance of momentum 
\begin{equation} 
\partial_t (\vr \vu) + \Div (\vr \vu \otimes \vu) + \Grad p(\vr) = 
\Div \mathbb{S}(\Ds \vu) + \vr \Grad G
\end{equation}
 that
\begin{equation} 
\partial_t \vr \vu+\vr\partial_t (\vu - \vu_E) + \Div (\vr \vu \otimes \vu) + \Grad p(\vr) = 
\Div \mathbb{S}(\Ds \vu) + \vr \Grad G
\end{equation}
and consequently
\begin{equation} \label{AP1}
\begin{split}
\partial_t \vr \vu\cdot (\vu - \vu_E)+\vr\partial_t (\vu - \vu_E)\cdot (\vu - \vu_E) + \Div (\vr \vu \otimes \vu)\cdot (\vu - \vu_E) + \Grad p(\vr)\cdot (\vu - \vu_E) &\\ 
= 
\Div \mathbb{S}(\Ds \vu)\cdot (\vu - \vu_E) + \vr \Grad G\cdot (\vu - \vu_E).
\end{split}
\end{equation}

Since
$\partial_t (\frac12 \vr |\vu - \vu_E|^2) = \frac12 \partial_t \vr |\vu - \vu_E|^2 + \vr (\vu - \vu_E)\cdot \partial_t (\vu - \vu_E) $ we rewrite \eqref{AP1} as
\begin{equation} 
\begin{split}
\frac12 \partial_t \vr |\vu - \vu_E|^2 - \frac12 \partial_t \vr |\vu - \vu_E|^2 +\partial_t \vr \vu\cdot (\vu - \vu_E)+\vr\partial_t (\vu - \vu_E)\cdot (\vu - \vu_E)& \\
+ \Div (\vr \vu \otimes \vu)\cdot (\vu - \vu_E) + \Grad p(\vr)\cdot (\vu - \vu_E) = 
\Div \mathbb{S}(\Ds \vu)\cdot (\vu - \vu_E) + \vr \Grad G\cdot (\vu - \vu_E)
\end{split}
\end{equation}
and then simplify it to
\begin{equation} \label{AP2}
\begin{split}
\partial_t (\frac12 \vr |\vu - \vu_E|^2) - \frac12 \partial_t \vr |\vu - \vu_E|^2 +\partial_t \vr \vu\cdot (\vu - \vu_E)& \\
+ \Div (\vr \vu \otimes \vu)\cdot (\vu - \vu_E) + \Grad p(\vr)\cdot (\vu - \vu_E) = 
\Div \mathbb{S}(\Ds \vu)\cdot (\vu - \vu_E) + \vr \Grad G\cdot (\vu - \vu_E).
\end{split}
\end{equation}

Next, we first integrate \eqref{AP2} over $\Omega$, then multiply by $\psi \in C^1_c[0, \infty)$, $\psi \geq 0$ and finally integrate over $(0,\infty)$ to get
\begin{equation} \label{AP3}
\begin{split}
\int_0^\infty \psi \int_\Omega \partial_t (\frac12 \vr |\vu - \vu_E|^2) \dx \dt - \int_0^\infty \psi \int_\Omega \frac12  \partial_t \vr |\vu - \vu_E|^2\dx \dt +\int_0^\infty \psi \int_\Omega\partial_t \vr \vu\cdot (\vu - \vu_E)\dx \dt& \\
+ \int_0^\infty \psi \int_\Omega\Div (\vr \vu \otimes \vu)\cdot (\vu - \vu_E)\dx \dt + \int_0^\infty \psi \int_\Omega\Grad p(\vr)\cdot (\vu - \vu_E) \dx \dt& \\
= \int_0^\infty \psi \int_\Omega \Div \mathbb{S}(\Ds \vu)\cdot (\vu - \vu_E)\dx \dt + \int_0^\infty \psi \int_\Omega \vr \Grad G\cdot (\vu - \vu_E)\dx \dt
\end{split}
\end{equation}

Using integration by parts we can rewrite the first term in \eqref{AP3} as
\begin{equation} \label{AP4}
\int_0^\infty \psi \int_\Omega \partial_t (\frac12 \vr |\vu - \vu_E|^2) \dx \dt = - \int_0^\infty \partial_t \psi \int_\Omega \frac12 \vr |\vu - \vu_E|^2 \dx \dt  - \psi(0) \int_\Omega \frac12 \vr(0,\cdot) |\vu(0,\cdot) - \vu_E|^2 \dx. 
\end{equation}

The rest of the terms in \eqref{AP3} can be treated as follows (below we ignore the time integration and multiplication by $\psi$ as it plays no role in the calculations):
\begin{itemize}
\item
\begin{equation}
\begin{split}
 \int_\Omega\Div (\vr \vu \otimes \vu)\cdot (\vu - \vu_E)\dx  &= -   \int_\Omega(\vr \vu \otimes \vu): \Grad(\vu - \vu_E)\dx  +  \int_{\partial \Omega} (\vr \vu \otimes \vu)\cdot (\vu - \vu_E)\cdot \vn dS  \\
&=  -   \int_\Omega(\vr \vu \otimes \vu): \Grad \vu \dx  +   \int_\Omega(\vr \vu \otimes \vu): \Grad  \vu_E\dx
\end{split}
\end{equation}
where the boundary term vanishes thanks to the boundary conditions  $\vu = \vuE$ on $\partial \Omega$.

\item

\begin{equation} 
\begin{split}
& -  \int_\Omega \frac12  \partial_t \vr |\vu - \vu_E|^2\dx  + \int_\Omega\partial_t \vr \vu\cdot (\vu - \vu_E)\dx  -   \int_\Omega(\vr \vu \otimes \vu): \Grad \vu \dx  \\
&=  -  \int_\Omega \frac12  \partial_t \vr (\vu - \vu_E)\cdot (\vu - \vu_E)\dx + \int_\Omega\partial_t \vr \vu\cdot (\vu - \vu_E)\dx  -   \int_\Omega(\vr \vu \otimes \vu): \Grad \vu \dx  \\
&=  \int_\Omega \frac12  \partial_t \vr \vu \cdot (\vu - \vu_E)\dx + \int_\Omega \frac12\partial_t \vr \vu_E\cdot (\vu - \vu_E)\dx  -   \int_\Omega(\vr \vu \otimes \vu): \Grad \vu \dx  \\
&= \frac12 \int_\Omega   \partial_t \vr  (\vu+ \vu_E) \cdot (\vu - \vu_E)\dx  -   \int_\Omega(\vr \vu \otimes \vu): \Grad \vu \dx  \\
&=  - \frac12 \int_\Omega   \Div(\vr \vu) (|\vu|^2 - |\vu_E|^2)\dx   -   \int_\Omega(\vr \vu \otimes \vu): \Grad \vu \dx  \\
&=   \frac12 \int_\Omega   \vr \vu \cdot (\Grad |\vu|^2 - \Grad|\vu_E|^2)\dx -  \frac12 \int_{\partial \Omega}   \vr \vu \cdot \vn ( |\vu|^2 - |\vu_E|^2)\dx  -   \int_\Omega(\vr \vu \otimes \vu): \Grad \vu \dx  \\
&=   - \frac12 \int_\Omega   \vr \vu \cdot \Grad|\vu_E|^2 \dx +  \frac12 \int_\Omega   \vr \vu\cdot \Grad |\vu|^2\dx  -   \int_\Omega(\vr \vu \otimes \vu): \Grad \vu \dx  \\
&=   - \frac12 \int_\Omega   \vr \vu \cdot \Grad|\vu_E|^2 \dx  \\
\end{split}
\end{equation}

\item 
\begin{equation}
\begin{split}
\int_\Omega \Div \mathbb{S}(\Ds \vu)\cdot (\vu - \vu_E)\dx &= - \int_\Omega  \mathbb{S}(\Ds \vu): \Grad (\vu - \vu_E)\dx + 
\int_{\partial \Omega} \mathbb{S}(\Ds \vu)\cdot (\vu - \vu_E)\cdot \vn dS\\
& = - \int_\Omega  \mathbb{S}(\Ds \vu): \Ds \vu \dx + \int_\Omega  \mathbb{S}(\Ds \vu): \Ds \vu_E\dx ,
\end{split}
\end{equation}
where the boundary term vanishes thanks to the boundary conditions on $\vu$ and 
\begin{equation}
\mathbb{S}(\Ds \vu): \Grad (\vu - \vu_E) = \mathbb{S}(\Ds \vu): \Ds (\vu - \vu_E)
\end{equation} thanks to the symmetry of $\mathbb{S}(\Ds \vu)$.

\item 

\begin{equation}
\begin{split} \label{AP5}
\int_\Omega\Grad p(\vr)\cdot (\vu - \vu_E) \dx &=\int_\Omega\Grad p(\vr)\cdot \vu \dx  + \int_\Omega p(\vr)\Div \vu_E \dx - \int_{\partial \Omega} p(\vr) \vu_E \cdot \vn dS \\
&= \int_\Omega\Grad p(\vr)\cdot \vu \dx  + \int_\Omega p(\vr) \mathbb{I}  : \Grad \vu_E \dx- \int_{\partial \Omega} p(\vr) \vu_E \cdot \vn dS.
\end{split}
\end{equation}

\end{itemize}

Finally, thanks to the boundary conditions on $\vu$ we get

\begin{equation}
\begin{split}
\int_\Omega \partial_t P(\vr) \dx & = \int_\Omega P'(\vr)\partial_t \vr \dx  = - \int_\Omega P'(\vr)\Div (\vr \vu)\dx\\
& = \int_\Omega \Grad P'(\vr) \cdot (\vr \vu) \dx - \int_{\partial \Omega} P'(\vr) \vr \vu \cdot \vn dS
 = \int_\Omega  P''(\vr) \vr \Grad \vr \cdot \vu \dx - \int_{\partial \Omega} P'(\vr) \vr \vu \cdot \vn dS\\
& =  \int_\Omega  p'(\vr) \Grad \vr \cdot \vu \dx - \int_{\partial \Omega} P'(\vr) \vr \vu \cdot \vn dS =  \int_\Omega \Grad p(\vr) \cdot \vu \dx - \int_{\partial \Omega} P'(\vr) \vr \vu_E \cdot \vn dS
\end{split}
\end{equation}
 and hence

\begin{equation}
\begin{split} \label{APL}
 - \int_0^\infty \partial_t \psi \int_\Omega P(\vr) \dx \dt &- \psi(0) \int_\Omega P(\vr(0,\cdot))\dx = \int_0^\infty \psi \partial_t \int_\Omega P(\vr) \dx \dt \\
& = \int_0^\infty \psi  \int_\Omega \Grad p(\vr) \cdot \vu \dx \dt  - \int_0^\infty \psi \int_{\partial \Omega} P'(\vr) \vr \vu_E \cdot \vn dS \dt,
\end{split}
\end{equation}
for any $\psi \in C^1_c[0, \infty)$, $\psi \geq 0$.

Relations \eqref{AP4}--\eqref{AP5}, \eqref{APL} put together with  \eqref{AP3} and the boundary conditions  $\vr=\vrE$ on $\Gamma_{\mbox{in}}$ yield the energy inequality \eqref{M3} for all $\psi \in C^1_c[0, \infty)$, $\psi \geq 0$.


\def\cprime{$'$} \def\ocirc#1{\ifmmode\setbox0=\hbox{$#1$}\dimen0=\ht0
  \advance\dimen0 by1pt\rlap{\hbox to\wd0{\hss\raise\dimen0
  \hbox{\hskip.2em$\scriptscriptstyle\circ$}\hss}}#1\else {\accent"17 #1}\fi}


\end{document}